\documentclass[11pt]{amsart}

\usepackage{amssymb}

\usepackage{xypic}

\newtheorem{theorem}{Theorem}[section]
\newtheorem{proposition}[theorem]{Proposition}
\newtheorem{lemma}[theorem]{Lemma}
\newtheorem{definition}[theorem]{Definition}

\newtheorem{corollary}[theorem]{Corollary}
\newtheorem{remark}[theorem]{Remark}

\numberwithin{equation}{section}

\begin{document}

\def\da{\!\downarrow\!}
\def\ua{\!\uparrow\!}

\def\ker{\hbox{ker }}
\def\NN{\mathbb N}
\def\ZZ{\mathbb Z}
\def\Ext{\hbox{Ext}}
\def\Hom{\hbox{Hom}}
\def\Ann{\hbox{Ann}}
\def\Rack{\hbox{Rack}}
\def\Fac{\hbox{Fac}}
\def\Ab{\hbox{Ab}}
\def\Gp{\hbox{Gp}}
\def\Mon{\hbox{Mon}}
\def\mod{\hbox{mod}}
\def\be{\begin{equation}}
\def\ee{\end{equation}}
\def\tensor{\otimes}
\def\iso{\cong}
\def\op{\mathrm{op}}
\def\Cat{\mathbf{Cat}}
\def\cM{\mathcal{M}}
\def\cB{\mathcal{B}}
\def\cGpd{\mathcal{G}}
\def\Set{\mathbf{Set}}
\def\ob{\mathrm{ob}\,}
\def\bB{\mathbf{B}}
\def\bM{\mathbf{M}}
\def\bC{\mathbf{C}}
\def\bG{\mathbf{G}}
\def\la#1{\mathop{\longleftarrow}\limits^{#1}}
\def\ra#1{\mathop{\longrightarrow}\limits^{#1}}
\def\colim{\mathop{\mathrm{colim}}}
\def\Cat{\mathbf{Cat}}
\def\pr{\mathrm{pr}}
\def\ff{\mathrm{ff}}
\def\inj{\mathrm{in}}
\def\Aut{\mathrm{Aut}}
\def\bA{\mathbf{A}}
\def\bSp{\mathbf{Sp}}
\def\Sp{\mathrm{Sp}}

\title{The Burnside bicategory of groupoids}

\author{Haynes Miller}
\address{Department of Mathematics\\
Massachusetts Institute of Technology\\
Cambridge, MA 02139}
\email{hrm@math.mit.edu}
\urladdr{http://math.mit.edu/$\sim$hrm}

\maketitle

\begin{center}
{\em Dedicated to the memory of Gaunce Lewis, 1950--2006,\\
and Sam Gitler, 1933--2014}
\end{center}

\medskip
Sometime around 1980 Gaunce Lewis described to me a ``Burnside category'' 
of groupoids. It provides an additive completion of the Burnside 
category of groups. In these notes we describe this
construction in the more structured setting of bicategories \cite{leinster}.

Morphisms from
$H$ to $G$ in the Burnside category of groups are provided by sets with 
commuting right and left actions of $H$ and $G$, which are free and finite
over $G$. We extend this construction to groupoids by
considering a $G$-action on a set $X$ as a homomorphism $G\rightarrow\Aut(X)$.
This has an alternative expression in terms of ``bi-actions,''  
a variant of
the discrete version of a topological construction which provides a convenient
account of the theory of orbifolds \cite{lerman}.
We obtain a bicategory $\bB$, the ``Burnside bicategory of groupoids.'' 
It comes equipped with a morphism \cite{leinster} from the bicategory $\bG$ of 
groupoids, functors, and natural transformations, and a compatible 
functor from the opposite of the sub-bicategory $\bG_{cf}$ generated by
``finite covers'' of groupoids. We think of the first as ``stabilization''
and the second as providing ``transfers.'' 

The morphism categories are 
symmetric monoidal, and the morphism
groups in the Burnside category (of groupoids, as of groups) 
are obtained by adjoining inverses
to the commutative monoids of isomorphism classes.

Moreover, the composition functors are ``bilinear'' and presumably
satisfy the axioms described by Bert Guillou \cite{guillou}. Guillou 
shows that such pairings can be rigidified to provide the input required by 
Tony Elmendorf and Mike Mandell \cite{elmendorf-mandell} 
to produce a spectral category $\Sp\bG$.
The functor $G\mapsto\Sigma^\infty_+BG$ will then
extend to a spectral functor to spectra. On the subcategory of finite
groupoids, the maps on morphism spectra
are appropriate completions; this is the content of the ``Segal conjecture.''

In fact the disjoint union of groupoids plays the role of both coproduct
and product in these bicategories. In this sense we have a ``semi-additive
bicategory.'' We have not attempted to write down 
a proper definition of this, but presumably: there is
such a notion (see \cite{tommasini} for some relevant work);
it implies Guillou's axioms (for a ``pre-additive bicategory'');
and examples such as the ones
presented here are additive bicategories. In any case, passing to the 
underlying Burnside category of groupoids we obtain an additive category with 
compatible functors from the category of groupoids and functors
and from the opposite of the subcategory in which morphisms are finite 
covering maps. 

This construction can be copied easily
to give a ``Morita category'' of rings, in which the morphism symmetric
monoidal categories are the groupoids of bimodules which are finitely generated
and projective over the source ring, and bimodule isomorphisms. 
This leads to an extension of the
construction of the algebraic K-theory of a ring to a spectral functor
from the ``spectral Morita category.'' 

An $H\times G$ bi-set is a functor $H^{\op}\times G\rightarrow\Set$. 
Such functors go by various names in the literature: ``modules,'' 
``bimodules,''
or ``profunctors.'' They are also sometimes called ``correspondences,'' but
it seems better to reserve that term for a ``span,'' a diagram of the form
$H\leftarrow K\rightarrow G$, with suitable conditions of the functor
$K\rightarrow H$. In the second part of this paper, we investigate a
bicategory $\bC$ of correspondences of groupoids. The appropriate condition
on $K\rightarrow H$ seems to be that it should be a ``weak finite cover''
(Definition \ref{def-wfc}). We set up this bicategory and show how it is
related to the bicategory $\bB$, by means of appropriate co-end and 
translation category constructions. The stabilization and transfer 
functors are tautologous in the setting of corespondences, and display
it as a bicategorical form of Quillen's $Q$-construction. 

My interest in making these constructions explicit was re-ignited by a lecture
by Clark Barwick in April, 2009, in which he carried out the analogue for
schemes, and showed that algebraic $K$-theory extends to an enriched
functor on the resulting spectral category. I claimed at the time
that the construction for groupoids was standard and well known, but 
was not able to back up that claim when Dustin Clausen asked what I meant. 
These notes are my attempt to make good on that claim. 
My objective is to be self-contained and explicit. Much deeper and far 
reaching work is being done in this area by many people, among them
Clark Barwick, Anna Marie Bohmann, David Gepner, Bert Guillou, Rune Haugseng,
Peter May, and Stefan Schwede on the topological side, 
and most recently Emily Riehl and Dominic Verity \cite{riehl-verity}
on the category theoretic side. 

The first section below describes various notions connected with actions
of groupoids on sets. This is extended in Section 2 to bi-sets, and is used
set up the ``Burnside bicategory'' of groupoids. We recall some standard
constructions in Section 3, and then describe the connection with topology
as well as the two connections with the usual bicategory of groupoids. 
Section 4 is devoted to the correspondence bicategory, and Section 5 to 
the relation between it and the bi-set bicategory.

I am grateful to Chris Schommer-Pries for pointing out the bi-action model
(see \cite{lerman,schommer-pries}).
Conversations with Jacob Lurie, Clark Barwick, Matthew Gelvin, 
Ang\'elica Osorno, and Martin Frankland have also been helpful,
and I thank the referee for a very careful reading and for pointing 
out a signal error in the first version of this paper, correcting
which led to an extensive reorganization of the paper.
My longtime friendships with Sam Gitler and with Gaunce Lewis
enriched my mathematical life over many years, and I am forever grateful.

\section{$G$-sets, $G$-actions, and covers}

We begin by establishing several equivalent ways to view an ``action''
of a groupoid on a set, and note how the conditions of being finite or
free appear from these various perspectives. To start with, we set up some
basic categorical features of the category of $G$-sets.

\subsection{$G$-sets}
Write $G_0$ for the set of objects of the groupoid $G$, and $G_1$ for
the set of morphisms.
Call a functor $X:G\rightarrow\Set$ from a groupoid $G$ to the category of
sets a ``left $G$-set'' (or just a ``$G$-set'')
and a functor $Y:G^\op\rightarrow\Set$ a ``right $G$-set.'' Write $X_\gamma$
and $Y^\gamma$ for the values at $\gamma\in G_0$. 
If $g:\gamma'\rightarrow\gamma$ is a morphism in $G$, write
$gx\in X_\gamma$ for the image of $x\in X_{\gamma'}$ under $\gamma$, and 
$yg\in Y^{\gamma'}$ for the image of $y\in Y^\gamma$. Left $G$-sets form
a category $G$-$\Set$ in which morphisms are natural transformations.

The forgetful functor from $G$-sets to $G_0$-sets (where we regard
$G_0$ as the discrete subgroupoid of $G$ consisting of the identity
morphisms) preserves limits and colimits. Moreover:

\begin{lemma} The forgetful functor $G$-$\Set\rightarrow G_0$-$\Set$
creates coequalizers.
\label{lemma-coequalizers} 
\end{lemma}

\noindent
{\em Proof.}
Let $u,v:X\rightarrow Y$, and for each $\gamma\in G_0$
let $Z_\gamma$ be the coequalizer of 
$u_\gamma,v_\gamma:X_\gamma\rightarrow Y_\gamma$. We need to see that 
$Z$ extends uniquely to a functor from $G$, so let 
$g:\gamma'\rightarrow\gamma$. Since $u$ and $v$ are natural
transformations, there is a unique map 
$g:Z_{\gamma'}\rightarrow Z_{\gamma}$ compatible with
$g:Y_{\gamma'}\rightarrow Y_{\gamma}$ under the projection map
$Y\rightarrow Z$. Uniqueness implies that with this structure $Z$ is
a $G$-set. To check that it is the coequalizer, let $w:Y\rightarrow W$
be a map of $G$-sets such that $wu=wv$. For each $\gamma$, there is a
unique map $Z_\gamma\rightarrow W_\gamma$ factorizing $f_\gamma$. 
These maps are compatible with $g$ because they are after composing
with the surjection $Y_{\gamma}\rightarrow Z_{\gamma}$. $\Box$

This lemma lets us recognize effective epimorphisms in $G$-$\Set$.
Recall that a map $X\rightarrow Y$ is an effective epimorphism if it is
the coequalizer of the canonical 
fork $X\times_YX\rightrightarrows X$. This is the
categorical notion of dividing by an equivalence relation.
The effective epimorphisms in $\Set$ are the surjections.

\begin{corollary}
A morphism in $G$-$\Set$ is an effective epimorphism if and only if it is
surjective on objects. 
\label{cor-effective-epi}
\end{corollary}

\noindent
{\em Proof.} Effective epimorphisms are preserved by the forgetful 
functor since it preserves limits and colimits. Conversely, 
suppose that $f:X\rightarrow Y$ is such that $f_\gamma$ is surjective
for each $\gamma\in G_0$. Surjections of sets are effective epimorphisms;
that is, 
\[
X_\gamma\times_{Y_\gamma}X_\gamma\rightrightarrows X_\gamma\rightarrow Y_\gamma
\]
is a coequalizer diagram. But the forgetful functor preserves pullbacks, so that
$X_\gamma\times_{Y_\gamma}X_\gamma=(X\times_YX)_\gamma$, and
by the lemma the diagram extends to a coequalizer diagram in $G$-$\Set$. $\Box$

\begin{definition}
A $G$-set $X$ is {\em free} 
if for all $\gamma,\gamma'\in G_0$ and all $x\in X_{\gamma'}$
the map $G(\gamma',\gamma)\rightarrow X_{\gamma}$ sending $g$ to $gx$ 
is injective. 
\label{def-free}
\end{definition}

For example,
the {\em empty groupoid} $\varnothing$ has no objects (and hence no morphisms).
It is initial in the category of categories. So there is a single 
$\varnothing$-set (namely the empty one), for which the 
freeness condition is satisfied vacuously. More interesting
examples are provided by the following lemma.

\begin{lemma} The following are equivalent conditions on a $G$-set $X$.\\
{\rm(1)} $X$ is a coproduct of co-representable $G$-sets.\\
{\rm(2)} $X$ is free.\\ 
{\rm(3)} $X$ is projective in the category $G$-$\Set$.
\label{lemma-free}
\end{lemma}

\noindent
{\em Proof.}
To begin with, notice that any $G$-set $X$ is canonically the target of an 
effective epimorphism from a coproduct of co-representable $G$-sets, 
given by the evident natural map 
\[
\coprod_{\gamma'\in G_0}\coprod_{x\in X_{\gamma'}}G(\gamma',-)\rightarrow X\,.
\]
The source of this map can be made smaller, however. Write 
\[
G\backslash X=\colim_GX\,.
\]
The natural map 
\[
\coprod_{\gamma'\in G_0}X_{\gamma'}\rightarrow G\backslash X
\]
is surjective, so by the axiom of choice it admits a section. 
This means that there are subsets $X'_{\gamma'}\subseteq X_{\gamma'}$
such that the composite map 
\[
\coprod_{\gamma'\in G_0}X'_{\gamma'}\rightarrow G\backslash X
\]
is bijective. The evident natural map
\be
i:\coprod_{\gamma'\in G_0}\coprod_{x\in X'_{\gamma'}}G(\gamma',-)\rightarrow X
\label{enough}
\ee
is then still an effective epimorphism. 

We aim to show that this map is an isomorphism if $X$ is free. 
First, let $\gamma',\gamma''\in G_0$, $x\in X'_{\gamma'}$,  
$y\in X'_{\gamma''}$, $f:\gamma'\rightarrow\gamma$ and 
$g:\gamma''\rightarrow\gamma$, and assume that 
$fx=gy$. Then $g^{-1}f:\gamma'\rightarrow\gamma''$ sends $x$ to $y$,
and so $x$ and $y$ have the same image in
$G\backslash X$. Therefore $\gamma'=\gamma''$ and $x=y$. 

If we now assume
that $X$ is free, this implies that $f=g$, and the map $i$ is injective 
as well as surjective on objects, and hence an isomorphism: so any free
$G$-set is a coproduct of co-representables. 

We will also need the converse: 
For any $\gamma_0\in G_0$, the co-representable functor $G(\gamma_0,-)$ is
free: Given $\gamma,\gamma'\in G_0$ and $x\in G(\gamma_0,\gamma')$,
$G(\gamma',\gamma)\rightarrow G(\gamma_0,\gamma)$ by $f\mapsto f\circ x$
is not just monic but bijective. So any coproduct of co-representable
$G$-sets is free. 

For any $\gamma_0\in G_0$, the corepresentable functor $G(\gamma_0,-)$ is
projective, and coproducts of corepresentables are too.

Finally, if $X$ is projective then the map in (\ref{enough}) is split-epi,
so $X$ is a subobject of a coproduct of co-representables and hence of 
a free object. But clearly any subobject of a free $G$-set is free. $\Box$

Two distinct finiteness conditions are relevant: 

\begin{definition} 
A $G$-set $X$ is {\em finitely generated} if its colimit  
is finite. It is {\em finite} if $X_\gamma$ is finite for all $\gamma\in G_0$.
\end{definition}

A $G$-set is finitely generated if and only if it is a quotient of a coproduct
of finitely many co-representable $G$-sets.

\subsection{Some constructions} 
We recall a few standard constructions on $G$-sets. 

The {\em co-end} of functors 
$X:G^{\op}\rightarrow\Set$ 
and $Y:G\rightarrow\Set$ is the set
\[
X\otimes_GY=\coprod_{\gamma\in G_0}X^\gamma\times Y_\gamma/\sim
\]
where the equivalence relation is given by $(xg,y)\sim(x,gy)$ with
$g:\gamma'\rightarrow\gamma$ in $G$, $x\in X^\gamma$, $y\in Y_{\gamma'}$.  

Note that $*\otimes_GY\ra{\iso}G\backslash Y$ 
(where $*$ denotes the constant functor with singleton values).

A functor $p:H\rightarrow G$ induces a pair of adjoint functors between the 
categories of $G$-sets and $H$-sets. The right adjoint sends a $G$-set 
$X:G\rightarrow\Set$ to the composite $Xp$. The left adjoint may be described
using the co-end: Given a right $H$-set $X$ or left $H$-set $Y$, there are
natural isomorphisms 
\[
(p_!Y)^\gamma\iso G(p,\gamma)\tensor_HY\,,\quad
(p_!X)_\gamma\iso X\tensor_HG(\gamma,p)\,.
\]
In particular,
\be
G(-,\gamma)\otimes_GY\iso Y_\gamma\,,\quad
X\otimes_GG(\gamma,-)\iso X^\gamma\,.
\label{unitors}
\ee
It's easy to check that
\[
p_!H(\eta,-)\iso G(p\eta,-)\,,\quad p_!H(-,\eta)=G(-,p\eta)\,.
\]
Combining these we have
\be
\label{adjoint1}
G(p,\gamma)\tensor_HH(\eta,-)\iso G(p\eta,\gamma)\,,\,
\ee
\be
\label{adjoint2}
H(-,\eta)\tensor_HG(\gamma,p)\iso G(\gamma,p\eta)\,.
\ee

\subsection{$G$-actions}
The notion of a $G$-set follows 
the image of a $G$-action on $X$ as a homomorphism $G\rightarrow\Aut(X)$.
Normally however one thinks of an action of a group $G$ on a set $X$ as a map
$\alpha:G\times X\rightarrow X$ satisfying certain properties. The
groupoid story can be developed following this model as well.

A {\em left action} of a groupoid $G=(G_0,G_1)$ is a set $P$ together
with a map $p:P\rightarrow G_0$ and an ``action map'' 
\[
\alpha:G_1{}^s\!\!\times_{G_0}P\rightarrow P
\]
over $G_0$ that is unital and associative. 
Here the prescript indicates that $G_1$ is to be regarded as a set over 
$G_0$ via the source map $s$. The fiber product is regarded as a set over 
$G_0$ via the target map $t$. 
In formulas, if we write $\alpha(g,x)=gx$, we are requiring
$p(gx)=tg, (g'g)x=g'(gx), 1x=x$. 
A morphism of left $G$-actions is defined in
the evident way.

A $G$-set $X:G\rightarrow\Set$ determines an action of $G$ on the set 
\[
P_X=\coprod_{\gamma\in G_0}X_\gamma
\]
by defining $gx$ to be $X_g(x)$. Conversely, an action of $G$
on $P$ determines a $G$-set $X_P$ by $(X_P)_\gamma=p^{-1}(\gamma)$, 
with functoriality determined in the evident way by the action map.
These constructions provide an equivalence of categories between $G$-sets
and $G$-actions. 

Under this correspondence, 
a $G$-set is finitely generated precisely when the orbit set $P/(x\sim gx)$ 
of the corresponding action is finite. The equivalence relation on $P$ is the
one given by $x\sim gx$. $P$ is finite when the fibers of $p:P\rightarrow G_0$
are finite. A $G$-set is free precisely when the corresponding
$G$-action satisfies the following property: the shear map
\[
\sigma:G_1{}^s\!\!\times_{G_0}P\rightarrow P\times P\quad,\quad
(g,x)\mapsto(x,gx)
\]
is injective. (Orbifold morphisms are modeled \cite{lerman} by bi-actions in 
which the shear map is required to be an isomorphism.)

\subsection{Covering maps of groupoids}
\label{ss-covering-maps}

\begin{definition} \cite{brown,higgins} 
A map of groupoids $p:H\rightarrow G$ is a 
{\em fibration}
provided that for any $\eta'\in H_0$ and any $g:p\eta'\rightarrow\gamma$ 
there exists $h:\eta'\rightarrow\eta$ in $H$ such that $ph=g$.
It is a {\em cover} provided there is only one such pair $(\eta,h)$,
and a {\em finite cover} provided that in addition the set $p^{-1}\gamma$ 
is finite for all $\gamma\in G_0$.
\end{definition}

\begin{remark}{\em
The notation $p^{-1}\gamma$ indicates the subcategory of the source with
objects mapping to $\gamma$ and morphisms mapping to $1_\gamma$; but it is
clear that if $p$ is a cover then this category has only identity 
morphisms, and nothing is lost by regarding it as a set. The ``fiber'' 
$p^{-1}\gamma$ is to be compared with the ``homotopy fiber'' $\gamma/p$. 
An object of $\gamma/p$ is a pair $(\eta\in H_0,g:\gamma\rightarrow p\eta)$.
A morphism $(\eta',g')\rightarrow(\eta,g)$ is a morphism 
$h:\eta'\rightarrow\eta$ such that $g=ph\circ g'$. There is a natural functor
$p^{-1}\gamma\rightarrow\gamma/p$ that sends $\eta$ to 
$(\eta,1_\gamma:\gamma\rightarrow p\eta)$. This functor is always full;
it is representative if $p$ is a fibration, and an equivalence if $p$ is
a cover. 
}
\label{homotopy-fiber}
\end{remark}

All three classes of maps are closed under composition and strict base-change. 

A $G$-set $X$ determines a new groupoid,
the {\em translation groupoid} $GX$ with 
\[
(GX)_0=\coprod_{\gamma\in G_0}X_\gamma
\] 
and $GX(x,y)=\{g\in G_1:gx=y\}$. The map of $G$-sets
$X\rightarrow*$ induces a functor $GX\rightarrow G$, which is a 
cover.

Conversely, a cover $p:H\rightarrow G$ determines a $G$-set 
$X$ with $X_\gamma=p^{-1}(\gamma)\subseteq H_0$ and 
$g:X_{\gamma'}\rightarrow X_\gamma$
defined using the unique morphism lifting. This establishes an equivalence
of categories between $G$-sets and covers of $G$. 

Our conditions
on $G$-sets translate to conditions on covers as follows. 
The orbit set or colimit $G\backslash X$ of a $G$-set is precisely the set of
components of the translation groupoid $GX$, so a cover $p:H\rightarrow G$ 
corresponds to a finitely generated $G$-set when the set of components of $H$ 
is finite. It corresponds to a finite $G$-set when it is a finite cover.
It corresponds to a free $G$-set when, for every
$\eta\in H_0$ and $\gamma\in G_0$, the map $G(p\eta,\gamma)\rightarrow H_0$
that sends $g$ to the target of its unique lift with source $\eta$ 
is injective.

Covers provide another source of free $H$-sets, 
as shown by the following easily checked lemma.

\begin{lemma} Let $p:H\rightarrow G$ be a cover of groupoids.
Then $G(\gamma,p)$ is a free $H$-set for all $\gamma\in G_0$, and if
$p$ is a finite cover then $G(\gamma,p)$ is finitely generated
as an $H$-set for all $\gamma\in G_0$.
\label{lemma-covering-free}
\end{lemma}

The second statement follows from the isomorphism
$H\backslash G(\gamma,p)\iso p^{-1}(\gamma)$. $\Box$

\section{Bi-sets and bi-actions}

Here we discuss the bivariant form of $G$-sets and $G$-actions, and show
how they form the objects of the morphism categories in a bicategory structure
on groupoids. 

\subsection{Bi-sets} For groupoids $G$ and $H$,
call a functor $H^\op\times G\rightarrow\Set$ an ``$(H,G)$-bi-set.''
An $(H,G)$-bi-set is {\em admissible} if $X^\eta$ is free and finitely 
generated as a $G$-set for every $\eta\in H_0$. 

Much of what we do would work just as well if we dropped the finiteness
condition; but it arises in the connection with topology, so we retain it 
anyway.

The admissible $(H,G)$-bi-sets form the objects of a category $\bB(H,G)$ in 
which morphisms are the natural {\em isomorphisms}. 
$\bB(H,G)$ is a symmetric monoidal groupoid, with the tensor product given by 
disjoint union.

This definition provides us with the morphism categories for a bicategory 
\cite{leinster}
in which the objects are groupoids. 
Define the composition functor $\bB(H,G)\times\bB(K,H)\rightarrow\bB(K,G)$ 
by sending $X,Y$ to $X\otimes_HY$ with
\[
(X\otimes_HY)^\kappa_\gamma=X_\gamma\otimes_HY^\kappa\,.
\]

We must check that $X\otimes_HY^\kappa$ is free and finitely generated. 
First freeness:
Since $X\otimes_H-$ commutes with coproducts, we may assume by Lemma
\ref{lemma-free} that $Y^\kappa=H(\eta,-)$ for some $\eta\in H_0$.
Then by (\ref{unitors}) $X\otimes_HY^\kappa\iso X^\kappa$, which is free by
assumption.

Now finite generation: Note that
\[
G\backslash(X\otimes_HY^\kappa)\iso*\otimes_G(X\otimes_HY^\kappa)\iso
(*\otimes_GX)\otimes_HY^\kappa\iso(G\backslash X)\otimes_HY^\kappa\,.
\]
Then use the facts that $G\backslash X^\eta$ is finite for every $\eta$,
that $H\backslash Y^\kappa$ is finite, and the following observation.

Let $X$ be a left $G$-set, and let $\{x_i\in X_{\gamma_i}\}$ 
represent the elements of $G\backslash X$. Then for any right $G$ set $W$ the
map 
\[
\coprod_i W^{\gamma_i}\rightarrow W\otimes_GX
\]
sending $w\in W^{\gamma_i}$ to $(w,\gamma_i,x_i)$ is surjective. 
To see this, let $(w,\gamma,x)\in W\times_GX$. 
There exists $i$ and $g:\gamma_i\rightarrow\gamma$ such that $gx_i=x$, so 
$(w,\gamma,x)=(w,\gamma,gx_i)=(wg,\gamma_i,x_i)$.

The identity object in $\bB(G,G)$ is given by 
\[
(1_G)^{\gamma'}_{\gamma}=G(\gamma',\gamma)
\]
The isomorphisms (\ref{unitors})  provide the unitors
\[
\rho:X\otimes_G1_G\rightarrow X
\quad,\quad
\lambda:1_G\otimes_GY\rightarrow Y
\]
The associator
\[
\alpha:X\otimes_G(Y\otimes_HZ)\rightarrow(X\otimes_GY)\otimes_HZ
\]
sends a list 
$(\varphi,x\in X_\varphi^\gamma,\gamma,y\in Y_\gamma^\eta,\eta,
z\in Z_\eta^\kappa,\kappa)$
to the same list, differently bracketed. The triangle identity for the
unitors and the pentagon identity for the associator are easily checked.

It is clear that 
\[
(X\amalg X')\tensor_GY\iso(X\tensor_GY)\amalg(X'\tensor_GY)\,,\,
X\tensor_G(Y\amalg Y')\iso(X\tensor_GY)\amalg(X\tensor_GY')\,.
\]
This is the main ingredient in the check that these pairings are
bilinear in the sense of \cite{guillou}.

\subsection{Bi-actions}
We now indicate the ``bi-action'' analogue of the bi-set bicategory. 

Given a right $G$-action on $Q$ and
a left $G$-action on $P$, the ``balanced product'' is defined as the set
\[
Q\times_GP=Q\times_{G_0}P/\sim
\]
where $(qg,p)\sim(q,gp)$ for $q\in Q^{\gamma}$,
$p\in P_{\gamma'}$, $g:\gamma'\rightarrow\gamma$. 

A {\em bi-action} of a pair of groupoids $G,H$ is a set
$P$ with a left action of $G$ (with projection $p:P\rightarrow G_0$) 
and a right action of $H$ (with projection $q:P\rightarrow H_0$), such
that $p(xh)=p(x)$, $q(gx)=q(x)$, and $(gx)h=g(xh)$. 

This is the same thing as a left action by $H^\op\times G$, so bi-actions
are equivalent to functors $H^\op\times G\rightarrow\Set$, i.e. to bi-sets.
The shear map for a bi-action is compatible with projections to $H_0$:
\[
\xymatrix{
G_1{}^s\!\!\times_{G_0}P \ar[rr]^\sigma \ar[dr] 
&& {P\times_{H_0}P} \ar[dl] \\
& {H_0} &
}
\]

Requiring $X^\eta$ to be free for every $\eta\in H_0$ is equivalent to
requiring that $\sigma$ be injective. Requiring $X^\eta$ to be finitely
generated for 
$\eta\in H_0$ is equivalent to requiring that in the corresponding 
bi-action $P$, the $G$-action on the fiber of $P$ over $\eta$ is finite.
Write $\bA(H,G)$ for the category of bi-actions of $(H,G)$ satisfying
these conditions. This is category is equivalent to $\bB(H,G)$.

Given another groupoid $K$, there is ``composition'' functor
$\bA(H,G)\times\bA(K,H)\rightarrow\bA(K,G)$ given by sending
$(P,Q)$ to the balanced product $Q\times_HP$, with the evident
actions of $K$ and $G$. We leave it to the reader to complete the structure
of a bicategory, and to verify that this bicategory is bi-equivalent
\cite{leinster} to the bicategory $\bB$. We will prefer the bi-set 
model.

\subsection{Additivity}
One of the most pleasing aspects of the bi-set bicategory $\bB$ is that it is
{\em additive}: the coproduct groupoid serves as both coproduct and product
in $\bB$. 

Let $G$ and $H$ be two groupoids.
Define the admissible $(G,G\amalg H)$-bi-set $X(1)$ by
\[
X(1)_{\gamma}=G(-,\gamma)\quad,\quad
X(1)_{\eta}=\varnothing\quad,\quad
\gamma\in G_0\,, \eta\in H_0\,.
\]
The bi-set $X(2)$ has a similar description. 
These morphisms induce an equivalence of categories
\[
\bB(G\amalg H,K)\rightarrow\bB(G,K)\times\bB(H,K)
\]
which establishes $G\amalg H$ as the coproduct in the bicategory $\bB$.

On the other hand, define the $(G\amalg H,G)$-bi-set $Y(1)$ by
\[
Y(1)^{\gamma}=G(\gamma,-)\quad,\quad
Y(1)^{\eta}=\varnothing\quad,\quad
\gamma\in G_0\,, \eta\in H_0\,,
\]
and define $Y(2)$ similarly. Since $G(\gamma,-)$ is co-representable, 
$Y(1)$ is admissible. Thus we have morphisms 
in the category $\bB$ that together induce a functor
\[
\bB(K,G\amalg H)\rightarrow\bB(K,G)\times\bB(K,H)
\]
This functor is also an equivalence
The groupoid $G\amalg H$, together with these structure morphisms,
serves as a product in the bicategory $\bB$.

\subsection{Burnside category} By passing to the commutative monoids of 
isomorphism classes of objects in the morphism symmetric monoidal categories
and then group-completing, we receive the {\em Burnside category} of groupoids.
It is preadditive by construction: abelian group structures on the morphism
sets are given, and composition is bilinear. Since group-completion commutes
with finite products, the constructions given above continue to provide 
products and coproducts in the Burnside category, which is therefore
actually additive.

\section{Topology}

The classifying space construction provides a functor from $\bG$ to the 
category of spaces. Under this functor, covers pass to covering spaces,
and finite covers to finite covering spaces. 

The topological interest in the Burnside bicategory arises from prospect
of using the bicategorical structure to form a spectrally enriched 
category of groupoids. As indicated in the Introduction, it is to be hoped
that the bicategorical structure described here produces a spectrally
enriched category $\Sp\bB$.

Given a spectrally enriched category $\bC$, an object $F$ in $\bC$ defines
a spectral functor $\bC\rightarrow\bSp$
to the category $\bSp$ of spectra sending $G$ to $\bC(F,G)$. Then we may define
\[ 
\bC(H,G)\rightarrow \bSp(\bC(F,H),\bC(F,G))
\]
as the adjoint of the
composition morphism $\bC(H,G)\wedge\bC(F,H)\rightarrow\bC(F,G)$.
In our example, taking $F$ to be the singleton groupoid yields a 
spectral functor $\Sp\bB\rightarrow\bSp$ sending $G$ to the spectrum
associated to the category (with coproducts) consisting of free 
and finitely generated $G$-modules. Since a finitely generated free 
$G$-set splits canonically into a coproduct of corepresentables, and
\[
\Aut(\coprod_i n_iG(\gamma_i,-))=\prod_i\Sigma_{n_i}\wr\Aut(\gamma_i)\,,
\]
we find 
\[
\Sp\bB(*,G)=\Sigma^\infty_+BG
\]
The resulting maps 
\[
\Sp\bB(H,G)\rightarrow F(\Sigma^\infty_+ BH,\Sigma^\infty_+ BG)
\]
are, by the Segal conjecture, appropriate completions when $G$ and $H$ are
finite. 

\subsection{Stabilization and transfer}
So you can think of $\bB$ as the stabilization of the bicategory of groupoids.
As such it should enjoy transfer maps as well. In this section we describe 
this construction.

Let $\bG$ denote the bicategory of groupoids, functors, and natural
transformations. (It is actually a 2-category.) Let $\bG_{fc}$ denote the
sub-bicategory with $\bG_{fc}(H,G)$ given by the groupoid of finite covers
$H\rightarrow G$ and natural transformations between them. 

There are morphisms \cite{leinster}  
(automatically ``homomorphisms,'' in our groupoid context) 
or pseudo-functors of bicategories
\[
S:\bG\rightarrow\bB\,,\quad T:\bG_{fc}^{\op}\rightarrow\bB 
\]
each of which is the identity on objects. The morphism $S$ sends
a functor $p:H\rightarrow G$ to the $(H,G)$-bi-set
$S(p)\in\bB(H,G)_0$ defined by
\[
S(p)_\gamma^\eta=G(p\eta,\gamma)\,,
\]
which is admissible since corepresentable functors are free and have 
singleton orbit sets,
while $T$ sends a finite cover $p:H\rightarrow G$ to the 
$(G,H)$-bi-set $T(p)\in\bB(G,H)_0$ defined by
\[
T(p)^\gamma_\eta=G(\gamma,p\eta)\,,
\]
which is admissible by Lemma \ref{lemma-covering-free}. 

On 2-cells, the morphism $S$ sends a natural transformation 
$\alpha:p'\rightarrow p$ to the isomorphism of $(H,G)$-bi-sets given by
$(\alpha^{-1}_\eta)^*:G(p'\eta,\gamma)\rightarrow G(p\eta,\gamma)$,
while $T$ sends $\alpha$ to 
$(\alpha_\eta)_*:G(\gamma,p'\eta)\rightarrow G(\gamma,p\eta)$.

Compatibility of $S$  with composition $K\ra{q}H\ra{p}G$ 
is given by the isomorphisms
\[
S(p)_\gamma\otimes_H S(q)^\kappa\ra{\iso}S(p\circ q)_\gamma^\kappa
\]
provided by \eqref{adjoint1}.
The required compatibilities with the associators (a hexagon) and 
unitors are easily verified. Similarly, compatibility of $T$ is given by
the isomorphisms
\[
T(q)_\kappa\tensor_HT(p)^\gamma\ra{\iso}T(p\circ q)^\gamma_\kappa
\]
provided by \eqref{adjoint2}.

These two morphisms are linked by a ``double coset formula'' of the 
following form. Let $p:H\rightarrow G$ and $q:K\rightarrow G$. Assume that
$q$ is a finite cover, and form the strict pullback diagram
\[
\xymatrix{
H\times_GK \ar[r]^{\overline p} \ar[d]^{\overline q} & K \ar[d]^q \\
H \ar[r]^p & G
}
\]
Then $\overline q$ is again a finite cover, and the pair $(q,p)$ defines a 
natural isomorphism of $(H,K)$-bisets
\[
\alpha:
S(\overline p)\tensor_{H\times_GK}T(\overline q)\ra{\iso}T(q)\tensor_GS(p)\,.
\]
For given $\kappa\in K_0$ and $\eta\in H_0$, this is an isomorphism
\[
K(\overline p,\kappa)\tensor_{H\times_GK}H(\eta,\overline q)
\ra{\iso}G(-,q\kappa)\tensor_GG(p\eta,-)\,. 
\]
Such a morphism is given by sending $(k,h)$ to $(qk,ph)$. Surjectivity follows
from the assumption that $q$ is a cover, and injectivity follows from
a careful inspection of the equivalence relations defining the co-ends.

These functors and natural transformations 
enjoy various further coherence properties. To describe
some of them, let's write $p_*=S(p)$ and $q^*=T(q)$. 

Suppose we have a diagram
\[
\xymatrix{
& {\cdot} \ar[r]^{\hat p} \ar[d]^{\overline q'} & {\cdot} \ar[d]^{q'} \\
{\cdot} \ar[r]^{\overline p'} \ar[d]^{\hat q} & 
{\cdot} \ar[r]^{\overline p} \ar[d]^{\overline q} &
{\cdot} \ar[d]^{q} \\
{\cdot} \ar[r]^{p'} & {\cdot} \ar[r]^{p} & {\cdot} 
}
\]
of pullback squares, and suppose that $q$ and $q'$ are finite covers. 
Then $\overline q$, $\overline q'$, and $\hat q$ are as well, and the
following pentagons of natural transformation are commutative. We leave the
composition compatibility morphisms, 
$(\overline p\overline p')_*\rightarrow \overline p_*\overline p'_*$ 
and so on, undenoted.

\[
\xymatrix@C=2pt{
(\overline p\overline p')_*\hat q^* \ar[rr]^\alpha \ar[d] &&
q^*(pp')_* \ar[d] &&&
\hat p_*(\overline q\overline q')^* \ar[rr]^\alpha \ar[d] && 
(qq')^*p_* \ar[d] \\
\overline p_*\overline p'_*\hat q^* \ar[dr]_{\alpha} && 
q^*p_*p'_* &&&
\hat p_*\overline q'{}^*\overline q^* \ar[dr]_\alpha && 
q'^*q^*p_* \\
& \overline p_*\overline q^*p'_* \ar[ur]_\alpha &&&&& 
q^*\overline p_*\overline q^* \ar[ur]_\alpha
}
\]

\subsection{Remark on the Mackey category}
\label{subsec-Mackey}
I owe the following remark to Rune Haugseng. 
It is actually more natural to consider, instead of the category $\bB(H,G)$, 
what one might call the {\em bivariant Mackey bicategory} $\bM(H,G)$. 
The objects are 
the same as the objects of $\bB(H,G)$, namely $(H,G)$-bi-sets $X$ such that
for each $\eta\in H_0$ the $G$-set $X^\eta$ is free and finitely generated.
For any pair $X,Y$, of objects, we 
define $\Hom(X,Y)$ to be the category of spans \cite{borceux}
over this pair. 
Thus an object in $\Hom(X,Y)$   
is a diagram $X\leftarrow Z\rightarrow Y$ in $\bB(H,G)$, and a morphism 
from $X\leftarrow Z'\rightarrow Y$ is an isomorphism $Z'\rightarrow Z$
that is compatible with the projections to $X$ and $Y$. 

By replacing $\Hom(H,G)$ with the set of isomorphism classes of objects in
$\Hom(H,G)$ one obtains what one might call the ``bivariant Mackey category'' 
$\overline{\bM}(H,G)$. 
We remark that if a map $f:X\rightarrow Y$ in $\bB(H,G)$ is such that the span
$X\la{=}X\ra{f}Y$ is an isomorphism in $\overline{\bM}(H,G)$, then $f$ is
itself an isomorphism. So $\bB(H,G)$ is the subcategory of isomorphisms in
$\overline{\bM}(H,G)$.

When $H$ is a finite group and $G$ is the trivial groupoid $1$, $\bM(H,1)$ 
is the usual Mackey category of finite right $H$-sets and spans. It has finite
coproducts, and so a corresponding spectrum, which is the fixed point
spectrum of the $H$-equivariant sphere spectrum. The category $\Hom(X,Y)$
also has coproducts, and hence a corresponding spectrum. These will
fit together to give a spectral enrichment of the bicategory $\bM(H,1)$,
and presheaves on that spectral category provide one model for $H$-equivariant
stable homotopy theory.

This makes it seem likely that
the various bicategories $\bM(H,G)$ assemble to a ``Burnside 
tricategory'' that will provide a model for 
``global equivariant stable homotopy theory''
in the sense of Bohmann or Schwede.

\section{Correspondences of Groupoids} 

The construction of compatible suspension and transfer homomorphisms for 
groupoids suggests an alternative construction of the target bicategory 
in terms of spans or correspondences. In this section we will propose a 
construction of a bicategory $\bC$ of correspondences and in the next 
section we will check that it is bi-equivalent to the bi-set bicategory
$\bB$. 

We begin with some well-known considerations
about the ``slice'' or ``overcategory'' in the bicategorical context.
This is quite general, but we will continue to speak of groupoids.

Suppose $G$ is a groupoid. An object of the 
over-bicategory $\bG/G$ is a map $p:L\rightarrow G$ in $\bG$. A morphism
from $p':L'\rightarrow G$ to $p:L\rightarrow G$ is given by a functor
$t:L'\rightarrow L$ together with a natural transformation 
$\theta:p'\rightarrow pt$. The composition
$(t'',\theta'')\circ(t',\theta')$ is given by 
$(t''\circ t',\theta''_{t'}\circ\theta')$. 
The identity morphism on $p:L\rightarrow G$ is $(1,1)$.

The morphisms from $p'$ to $p$ form the objects of a category
$\bG/G(p',p)$, in which a morphism 
$(\overline t,\overline\theta)\rightarrow(t,\theta)$ consists in a 
natural transformation $\psi:\overline t\rightarrow t$ such that 
$p\psi\circ\overline\theta=\theta:p'\rightarrow pt$.
Composition is given by composition of natural transformations.

The composition law $(t'',\theta'')\circ(t',\theta')$ extends to a
functor $\bG/G(p,p'')\times\bG/G(p',p)\rightarrow\bG/G(p',p'')$.
Let $\psi':(\overline t',\overline\theta')\rightarrow(t',\theta')$,
$\psi'':(\overline t'',\overline\theta'')\rightarrow(t'',\theta'')$ be
morphisms. The composed morphism 
$\psi:
(\overline t''\overline t',\overline\theta''_{\overline t'}\overline\theta')
\rightarrow(t''t',\theta''_{t'}\theta')$
is given by the diagonal in the commutative diagram
\[
\xymatrix@=30pt{
\overline t''\overline t' 
\ar[r]^{\overline t''\psi'} \ar[d]_{\psi''_{\overline t'}} \ar[dr]^\psi &
\overline t''t' \ar[d]^{\psi''_{t'}} \\
t''\overline t' \ar[r]^{t''\psi'} & t''t'
}
\]

The required diagram commutes by virtue of the commutativity of  
\[
\xymatrix{
p' \ar[r]^{\overline\theta'} \ar[dd]_{\theta'} & 
p\overline t' \ar[r]^{\overline\theta''_{\overline t'}} \ar[ddl]^{p\psi'} &
p''\overline t''\overline t' \ar[d]^{p''\overline t''\psi'} \\
&& p''\overline t''t' \ar[d]^{p''\psi''_{t'}} \\
pt' \ar[rru]^{\overline\theta''_{t'}} \ar[rr]^{\theta''_{t'}} &&
p''t''t'
}
\]

\begin{definition}
\label{def-wfc}
A groupoid is {\em weakly discrete} if there is at most one 
morphism between any two objects---that is, each component is ``unicursal.'' 
A weakly discrete groupoid is {\em weakly finite} if it has finitely many 
components. A map $p:H\rightarrow G$ is a {\em weak cover} if it
is a fibration and for every $\gamma\in G_0$ the fiber 
$p^{-1}\gamma$ is a weakly discrete groupoid, and a {\em weak finite cover}
if all the fibers are weakly finite as well.
\end{definition}

Notice that the condition of being weakly discrete or weakly discrete and 
weakly finite is invariant under equivalence of groupoids. 
Weak covers and weak finite covers pull back to the same.
It's also easy to see that 
compositions of weak finite covers are weak finite covers, and that if 
$K'\rightarrow G$ and $K''\rightarrow G$ are
both weak finite covers then so is $K'\amalg K''\rightarrow G$.

We are now in a position to define the groupoid of correspondences
between two groupoids. I am grateful to the referee for pointing out
that my first attempt at this was not matching $\bB(H,G)$ but was
closer to the tricategorical speculation described in Section 
\ref{subsec-Mackey} above.

A ``correspondence'' from $H$ to $G$
is a diagram $H\la{q}L\ra{p}G$ of groupoids in which
$q$ is a weak finite cover. A morphism from
$(q',L',p')$ to $(q,L,p)$ is an equivalence class of triples 
$(\phi,t,\theta)$, 
where $t:L'\rightarrow L$ is a functor
of groupoids, $(t,\phi):L'\rightarrow L$ in the groupoid $\bG/H$,
and $(t,\theta):L'\rightarrow L$ in the groupoid $\bG/G$. Two triples
$(\overline\phi,\overline t,\overline\theta)$ and $(\phi,t,\theta)$ 
are equivalent if there is a natural 
transformation $\psi:\overline t\rightarrow t$ such that 
$q\psi\circ\overline\phi=\phi$ and $p\psi\circ\overline\theta=\theta$. 
We leave it to the reader to verify that the composition operation defined 
in $\bG/G$ passes to equivalence classes and defines a category. 

\begin{definition}
The {\em correspondence groupoid} $\bC(H,G)$ is the subcategory of isomorphisms
in the category of correspondences from $H$ to $G$. 
\end{definition}

The equivalence class of 
a morphism $(\phi,t,\theta):(q',L',p')\rightarrow(q,L,p)$ 
is an isomorphism exactly when $(\phi,t,\theta)$
is an ``internal equivalence'' \cite{leinster} in the bicategory 
$\bG/(H\times G)$; that is to say, when there exists 
$(\phi',t',\theta'):(q,L,p)\rightarrow(q',L',p')$ 
and natural transformations $\psi:1\rightarrow tt':L\rightarrow L$ and 
$\psi':1\rightarrow t't:L'\rightarrow L'$
such that each of the following diagrams commutes.
\be
\xymatrix{
p' \ar[r]^{\theta} \ar[dr]_{p'\psi'} & pt \ar[d]^{\theta'_t} &&
q' \ar[r]^{\phi} \ar[dr]_{q'\psi'} & qt \ar[d]^{\phi'_t} \\
& p't't &&& q't't\\
p \ar[r]^{\theta'} \ar[dr]_{p\psi} & p't' \ar[d]^{\theta_{t'}} &&
q \ar[r]^{\phi'} \ar[dr]_{q\psi} & q't' \ar[d]^{\phi_{t'}} \\
& ptt' &&& qtt'
}
\label{compatibility}
\ee

These groupoids form the morphism categories in a bicategory $\bC$,
the {\em correspondence bicategory} of groupoids. The composition functor
$\bC(G,F)\times\bC(H,G)\rightarrow\bC(H,F)$ is defined
using the strict pullback: 
\[
\xymatrix @! =15pt  {
&& X\times_GL \ar[dl]_{\overline n} \ar[dr]^{\overline p} \\
& L \ar[dl]_q \ar[dr]^p && K \ar[dl]_n \ar[dr]^m \\
K && G && F
}
\]
The pulled back functor $\overline n$ is a weak finite cover because
$n$ is, and so the composite $q\circ\overline n$ is again a weak finite
cover. 

We leave a discussion of the associator and unitors to the interested reader.

Stabilization and the transfer construction are even more tautologous
in this setting than when expressed in terms of bimodules.
The bicategory morphism $S:\bG\rightarrow\bC$ (for ``stabilization'') 
that  sends a groupoid 
to itself, a functor $p:H\rightarrow G$ to the object 
$H\la{1}H\ra{p}G$ of $\bC(H,G)$, and a natural transformation
$\theta:p'\rightarrow p$ to the class of $(1,1,\theta)$. 
The transfer construction is a morphism $T:\bG_{wf}^{\op}\rightarrow\bC$
from the full sub bicategory of weak finite covers. It sends a groupoid 
to itself, a weak finite cover $q:H\rightarrow G$ to the object 
$G\la{q}H\ra{1}H$ of $\bC(G,H)$, and a natural transformation 
$\phi:q'\rightarrow q$ to the equivalence class of $(\phi,1,1)$. 

These two functors are related by a ``double coset formula'':
Let $p:H\rightarrow G$ be a functor and $q:F\rightarrow G$ a weak finite
cover, and form the pullback $H\la{\overline q}F\times_GH\ra{\overline p}F$.
Write $S(p)$ and $T(q)$ for the apexes of the spans, so that we have two 
composed spans:
\[
\xymatrix @! =14pt {
&& T(q)\times_G S(p) \ar[dl] \ar[dr] 
&&&&& S(\overline p)\times_{F\times_GH}T(\overline q) \ar[dl] \ar[dr] \\
& S(p) \ar[dl] \ar[dr] && T(q) \ar[dl] \ar[dr] 
&&& T(\overline q) \ar[dl] \ar[dr] && S(\overline p) \ar[dl] \ar[dr] \\
H && G && F & H && F\times_GH && F 
}
\]

It is easy to construct a natural isomorphism
\[
T(q)\times_GS(p)\iso S(\overline p)\times_{F\times_GH}T(\overline q)
\]
in the groupoid $\bC(H,F)$. 
The pair of functors $S,T$, together with this 
natural isomorphism, is in a certain sense  universal among such. 
The bicategory $\bC$ is a bicategorical form of Quillen's $Q$-construction.

\section{Comparison}

Now we will see how the bicategories $\bB$ and $\bC$ are related.

It may help to look first at the case of groups.
So let $G$ be a group and $K$ 
a subgroup. The inclusion $p:K\rightarrow G$ expresses
$K$ as a groupoid over $G$; but there is another groupoid over $G$ that
is expressed in terms of the transitive $G$ set $G/K$, namely
the translation groupoid $G(G/K)$. These two groupoids over $G$ are
in a suitable sense equivalent. 

There is a natural map $\alpha:K\rightarrow G(G/K)$ in $\bG/G$ 
that sends the unique
object in $K$ to $K\in(G(G/K))_0$. Note that in $G(G/K)$, 
$\Aut(K)=\{g\in G:gK=K\}=K$, and this isomorphism 
defines the functor on morphisms.
This functor commutes strictly with the projection maps to $G$.

The map $\alpha$ has a non-natural quasi-inverse $\beta$ obtained 
by picking coset 
representatives $g_i$, so that $G=\coprod_i g_iK$. In terms of these, 
define $\beta:G(G/K)\rightarrow K$ by sending $g_iK$ to the unique object
of $K$, and $g:g_iK\rightarrow g_jK$ to $g_j^{-1}gg_i\in K$. This is 
functorial, but it does not commute with the projections to $G$.
However, there is a natural transformation 
$\theta:q\rightarrow p\beta$ (where $q:G(G/K)\rightarrow G$ is the projection
which forgets the objects). It is defined by $\theta_{g_iK}=g_i^{-1}$.

To simplify the formulae, we might as well chose the identity element of $G$ 
as the representative of the identity coset $K$. With this choice, the
composite $K\ra{\alpha}G(G/K)\ra{\beta}K$ is the identity morphism.

The composite $G(G/K)\ra{\beta}K\ra{\alpha}G(G/K)$ sends every object
$g_iK$ to $K$, and $g:g_iK\rightarrow g_jK$ to $g_j^{-1}gg_i$.
The natural transformation $q\rightarrow q\alpha\beta$ is given by 
$\theta_{g_iK}=g_i^{-1}$. 

This morphism in $\bG/G$ is far from an isomorphism, but there is an 
isomorphism from the identity to it in the category of endomorphisms
of $G(G/K)\rightarrow G$ in the bicategory $\bG/G$, given by 
$\psi_{g_iK}=g_i$.

Now we describe analogous constructions extending in two directions: to 
groupoids and to the bivariant situation. 

Let $G$ and $H$ be groupoids and $X$ an $(H,G)$-bi-set.
The ``double translation groupoid'' $GXH$ has objects 
\[
(GXH)_0=\coprod_{(\gamma,\eta)\in G_0\times H_0}X_\gamma^\eta
\]
and morphisms 
$(\overline\gamma,\overline x,\overline\eta)\rightarrow(\gamma,x,\eta)$ 
given by
pairs 
$g:\overline\gamma\rightarrow\gamma,h:\overline\eta\rightarrow\eta$ such that
$g\overline x=xh\in X_{\gamma}^{\overline\eta}$. 
Composition is evident: Given also
$(g',h'):(\gamma,x,\eta)\rightarrow(\gamma',x',\eta')$, 
the composite is given by $(g'g,h'h)$.

The double translation groupoid comes equipped with functors 
$p:GXH\rightarrow G$ and $q:GXH\rightarrow H$, both covariant. 

To better understand the double translation category, 
it is useful to re-express it in terms of the Grothendieck construction.
Let $F:H^{\op}\rightarrow\Cat$ be a functor. 
The {\em Grothendieck construction} of $F$ is the category $FH$ with objects
\[
(FH)_0=\coprod_{\eta\in H_0} F^\eta_0
\]
and morphisms 
$(\eta',\varphi'\in F^{\eta'}_0)\rightarrow(\eta,\varphi\in F^{\eta}_0)$
given by $(h:\eta'\rightarrow\eta,f:\varphi'\rightarrow\varphi h)$. 
There is a covariant  functor $FH\rightarrow H$ sending $(\eta,\varphi)$ to 
$\eta$. 
If $F^\eta$ has only identity morphisms, for each $\eta$, then the
Grothendieck construction coincides with the translation category. 
The evident projection functor $FH\rightarrow H$ is a always fibration.

Similarly, there is a Grothendieck construction $GF$ for a functor
$F:G\rightarrow\Cat$. If $X$ is an $(H,G)$-biset, then it is easy to construct
natural isomorphisms
\[
(GX)H\iso GXH \iso G(XH)
\]
Since the fiber over $\eta$ of $GXH\rightarrow H$ is the translation groupoid 
$GX^\eta$, we see:

\begin{lemma} If $X^\eta$ is free and finitely generated as a $G$-set for all 
$\eta\in H_0$ then $q:GXH\rightarrow H$ is a weak finite cover.  $\Box$
\end{lemma}

So far we have verified that the double translation groupoid construction 
gives a map on objects from $\bB(H,G)$ to $\bC(H,G)$. This construction 
extends to a functor: Let
$f:X'\rightarrow X$ be an isomorphism of $(H,G)$-bi-sets. It induces a functor
$t:GX'H\rightarrow GXH$ sending the object 
$(\gamma,x',\eta)$ to $(\gamma,fx',\eta)$ and the morphism
$(g,h):(\overline\gamma,\overline x',\overline\eta)\rightarrow(\gamma,x',\eta)$
to the morphism 
$(g,h):(\overline\gamma,f\overline x',\overline\eta)
\rightarrow(\gamma,fx',\eta)$: 
$g(f\overline x')=f(g\overline x')=f(x'h)=(fx')h$:
The natural transformations needed to define a morphism in 
$\bC(H,G)$ from $GX'H$ to $GXH$ can both be taken to be identity 
transformations.

We have now constructed a functor $\bB(H,G)\rightarrow\bC(H,G)$ for each
pair of groupoids $H,G$. This is the start of the description of a morphism
of bicategories \cite{leinster}, and the main result of this section is this:

\begin{theorem}
\label{relate}
These functors participate in a morphism of bicategories $\bB\rightarrow\bC$
that is in fact a bi-equivalence.
\end{theorem}

The construction of this bi-equivalence occupies the rest of this section.
By \cite{leinster}, it will suffice to construct a morphism that induces
equivalences of morphism categories and is surjective up to equivalence on
objects. The latter condition is obvious! 

To specify a morphism of bicategories one needs to specify a transformation
natural in the triple $G,H,K$
\be
\label{homomorphism}
[s]:(GXH)\times_H(HYK)\ra{} G(X\otimes_HY)K
\ee
in $\bC(H,K)$ (as well as a natural transformation comparing units which we
leave to the reader). 
The map $[s]$ is represented by a functor $s$ given on objects by sending
\[
((\gamma,x\in X_\gamma^\eta,\eta),
(\eta,y\in Y_\eta^\kappa,\kappa))
\]
to
\[
(\gamma,[x,y]\in(X\tensor_HY)^\kappa_\gamma,\kappa)\,.
\]
A morphism 
\[
((\overline\gamma,\overline x,\overline\eta),
(\overline\eta,\overline y,\overline\kappa))\rightarrow
((\gamma,x,\eta),(\eta,y,\kappa))
\]
in the pullback $(GXH)\times_H(HYK)$ consists of morphisms
$g:\overline\gamma\rightarrow\gamma$, 
$h:\overline\eta\rightarrow\eta$, $k:\overline\kappa\rightarrow\kappa$
such that $g\overline x=xh$ and $h\overline y=yk$.
This triple determines a morphism 
\[
(\overline\gamma,\overline x,\overline y,\overline\kappa)
\rightarrow(\gamma,x,y,\kappa)
\]
in $G(X\otimes_HY)K$ given by $(g,k)$. The required equality 
$(g\overline x,\overline y)=(x,yk)$ in $X\otimes_HY$
is established by the morphism $h$.

This functor is strictly compatible with the projections to $K$ and $G$.
It is not bijective on objects, but it does give an isomorphism of 
correspondences and hence a map in the
category $\bC(K,G)$. To construct an inverse we begin by constructing
a quasi-inverse $s'$ for the functor $s$. Begin by choosing,
for each $(\kappa,\gamma)\in K_0\times G_0$, a splitting of the surjection
\[
\coprod_{\eta\in H_0}X_\gamma^\eta\times Y_\eta^\kappa\rightarrow
(X\tensor_HY)^\kappa_\gamma
\]
---so for each class $[x,y]\in(X\tensor_HY)^\kappa_\gamma$, pick $\eta\in H_0$
and a representative $(x,y)\in X_\gamma^\eta\times Y_\eta^\kappa$. This
defines the quasi-inverse $s'$ on objects.
We need to define the quasi-inverse on the morphism
$g:\overline\gamma\rightarrow\gamma, k:\overline\kappa\rightarrow\kappa$
with $(g,k):[\overline x,\overline y]\mapsto[x,y]$ (where $x\in X_\gamma^\eta$,
$y\in Y_\eta^\kappa$, $\overline x\in X_{\overline\gamma}^{\overline\eta}$,
$\overline y\in Y_{\overline\eta}^{\overline\kappa}$ are the chosen 
representatives). The fact that $(g,k)$ is a morphism provides us with
$h:\overline\eta\rightarrow\eta$ such that $g\overline x=xh$ and
$h\overline y=yk$. Since $X_\gamma$ is a free right $H$-set, the pair
$(\overline\eta,h)$ is actually uniquely determined by the first equation.
Declare $s'(g,k)=(g,h,k)$.

It is clear that $ss'=1$; we can take $\psi=1$ in the definition of equivalence
of correspondences. The natural transformation $\psi':1\rightarrow s's$ 
is given on the object $(x,y)\in X_\gamma^\eta\times Y_\eta^\kappa$ by
$(1,h^{-1},1)$, where $h:\overline\eta\rightarrow\eta$ is the morphism 
described above relating the pair $(x,y)$ to the chosen representative
$s's(x,y)=(\overline x,\overline y)$ of $[x,y]\in X_\gamma\tensor_H Y^\kappa$.
The naturality diagram is assured again by the fact that $X_\gamma$ is
free as an $H$-set.

We leave to the reader the verification of the required compatibilities of the
natural transformations $[s]$ and the unital analogue.

Up to this point we have been able to use functors between the summits of
spans which strictly commute with the projection maps. 
The next step is to construct a quasi-inverse of the functor 
$\bB(H,G)\rightarrow\bC(H,G)$, and this will require the weaker notion
of morphism. 
In the case of groups, with $G=1$
you want to send a finite index subgroup $K$ of $H$ to the $H$-set $H/K$. 
When $G\neq1$ we get the $(H,G)$-set $H\times_KG$. 
The groupoid version sends $H\la{q}K\ra{p}G$ to the $(H,G)$-set $X$ with
$X_\gamma^\eta$ defined as the coend
\be
X_\gamma^\eta=G(p,\gamma)\otimes_KH(\eta,q)\,.
\label{b-morphism}
\ee

\begin{lemma}
If $q$ is a weak finite cover then $X^\eta$ is free and
finitely generated for all $\eta\in H_0$. 
\end{lemma}

\noindent
{\em Proof.}
To see  that $X^\eta$ is free we must show that for all
$x=(g:p\kappa\rightarrow\gamma,\kappa,h:\eta\rightarrow q\kappa)
\in X^\eta_\gamma$,
the map $G(\gamma,\gamma')\rightarrow X^\eta_{\gamma'}$ sending
$\overline g$ to 
$\overline gx=(\overline gg,\kappa,h)$ is injective. 
Suppose that $\overline g'x=\overline g x\in X^\eta_{\gamma'}$.
By the relation in the definition of the coend, there exists 
$k:\kappa\rightarrow\kappa$ such that $q(k)h=h$ and
$\overline g'gp(k)=\overline gg$. The first identity implies that
$k=1$ since by hypothesis $\eta/q$ is weakly discrete.
(See Remark \ref{homotopy-fiber}.)
Then the second relation implies that $\overline g=\overline g'$.

To see that $X^\eta$ is finitely generated note that
$*\otimes_GG(\gamma_0,-)=\pi_0(GG(\gamma_0,-))$  is a singleton 
because the translation category is unicursal, so (using an evident notation)
\[
*\otimes_{\gamma\in G}X^\eta_\gamma
=(*\otimes_{\gamma\in G}(G(p,\gamma))\otimes_KH(\eta,q)
=*\otimes_KH(\eta,q)\,.
\]
Now $KH(\eta,q)$ is the homotopy fiber of $q$ over $\eta$, and its set
of components $*\otimes_KH(\eta,q)$ is finite since $q$ is a weak finite  
cover. $\Box$

This construction extends to a functor $\bC(H,G)\rightarrow\bB(H,G)$. Let 
$(t:K'\rightarrow K$, $\theta:p'\rightarrow pt$, $\phi:q'\rightarrow qt)$
be a morphism
$(q',K',p')\rightarrow(q,K,p)$. To define the induced map
$G(p',\gamma)\otimes_{K'}H(\eta,q')\rightarrow
G(p,\gamma)\otimes_KH(\eta,q)$, pick $\kappa'\in K'_0$ and look at
the natural map $\inj_{\kappa'}:G(p'\kappa',\gamma)\times H(\eta,q'\kappa')
\rightarrow G(p',\gamma)\otimes_{K'}H(\eta,q')$. The map we are looking
for is induced by 
$\inj_{t\kappa'}\circ((\theta_{\kappa'}^{-1})^*\times(\phi_{\kappa'})_*):
G(p'\kappa',\gamma)\times H(\eta,q'\kappa')\rightarrow
G(p,\gamma)\otimes_KH(\eta,q)$. It is straightforward to check that 
the equivalence relation defining morphisms in $\bC(G,H)$ perfectly matches
the equivalence relation defining elements in the coend \eqref{b-morphism}.

\begin{proposition}
These two constructions are the functors in an adjoint equivalence
between $\bB(H,G)$ and $\bC(H,G)$. 
\end{proposition}

\noindent
{\em Proof.} Given an $(H,G)$-bi-set $X$, there is a map
\[
\beta:G(p,\gamma)\otimes_{GXH}H(\eta,q)\rightarrow X_\gamma^\eta
\]
given as follows. An object of $GXH$ is given as 
$(\overline\gamma,x\in X_{\overline\gamma}^{\overline\eta},\overline\eta)$. 
The map is
\[
(g:\overline\gamma\rightarrow\gamma,(\overline\gamma,x,\overline\eta),
h:\overline\eta\rightarrow\eta)\mapsto gxh\in X_\gamma^\eta
\]
This map is an isomorphism. It is surjective since
$(1,(\gamma,x,\eta),1)\mapsto x\in X^\eta_\gamma$.
To check that it is injective, let 
$(g':\overline\gamma'\rightarrow\gamma,(\overline\gamma',x',\overline\eta'),
h':\eta\rightarrow\overline\eta')$ map to the same element of $X^\eta_\gamma$:
$g'x'h'=gxh$. We seek a morphism in 
$GXH$---$(\hat g:\overline\gamma'\rightarrow\overline\gamma,
\hat h:\overline\eta'\rightarrow\overline\eta)$ with 
$\hat g x'=x\hat h\in X_{\overline\gamma}^{\overline\eta'}$---such that
$g\hat g=g'$ and $\hat hh'=h$. There is a unique solution, 
$\hat g=g^{-1}g',\hat h=hh'^{-1}$. Then we calculate that
$\hat gx'=g^{-1}g'x'=g^{-1}gxhh'^{-1}=x\hat h$.

Given a correspondence $H\la{q}K\ra{p}G$, there is a functor
\[
\alpha:K\rightarrow G(G(p,-)\otimes_KH(-,q))H
\]
defined as follows. An object in the double translation category is determined 
by $(\gamma\in G_0,g:p\kappa\rightarrow\gamma,\kappa\in K_0,
h:\eta\rightarrow q\kappa,\eta\in H_0)$. 
We send $\kappa\in K_0$ to
$(p\kappa,1,\kappa,1,q\kappa)$.
On morphisms, $k:\overline\kappa\rightarrow\kappa$ is sent to the pair
$(pk,qk)$; the coherence morphism is provided by $k$ itself.
This functor $\alpha$ is a morphism in the category of groupoids over
$G$ and over $H$; in fact the natural transformation allowed in a morphism
is the identity map in both cases. 

The required compatibilities are easily checked. 

The functor $\alpha$ is not an isomorphism of groupoids, 
but it is an isomorphism in $\bC(H,G)$. To see this, observe first that it
factors as
\[
K\ra{t}GG(p,-)K\times_KKH(-,q)H\ra{s}G(G(p,-)\tensor_KH(-,q))H
\]
where $s$ is the equivalence of (\ref{homomorphism}) and $t$ is the
functor sending $\kappa$ to 
$(1:p\kappa\rightarrow p\kappa,\kappa,1:q\kappa\rightarrow q\kappa)$
and a morphism $k:\overline\kappa\rightarrow\kappa$ to $(pk,k,qk)$. 
We construct a quasi-inverse $t'$ for $t$ that is compatible with the
projection functors 
\[
\xymatrix@C=10pt{
&GG(p,-)K\times_KKH(-,q)H \ar[dl]_{\hat q} \ar[dr]^{\hat p}\\
H && G
}\,.
\]

The functor $t'$ sends 
$(\gamma,g:p\kappa\rightarrow\gamma,\kappa,h:\eta\rightarrow q\kappa,\eta)$
to $\kappa$, and a morphism 
\[
(g',k',h'):(\overline\gamma,
\overline g:p\overline\kappa\rightarrow\overline\gamma,\overline\kappa,
\overline h:\overline\eta\rightarrow q\overline\kappa,\overline\eta)\rightarrow
(\gamma,g:p\kappa\rightarrow\gamma,\kappa,h:\eta\rightarrow q\kappa,\eta)
\]
(where $g':\overline\gamma\rightarrow\gamma$, 
$k':\overline\kappa\rightarrow\kappa$, and $h':\overline\eta\rightarrow\eta$
make
\[
\xymatrix{
p\overline\kappa \ar[r]^{\overline g} \ar[d]^{pk'} & 
\overline\gamma \ar[d]^{g'} &&
\overline\eta \ar[r]^{\overline h} \ar[d]^{h'} &
q\overline\kappa \ar[d]^{qk'} \\
pk \ar[r]^g & \gamma && \eta \ar[r]^{h} & q\kappa
}
\]
commute) to $k'$. This functor does not satisfy $qt'=\hat q$ or
$pt'=\hat p$, but there are natural transformations 
\[
\theta':\hat p\rightarrow pt'\,,\quad\phi':\hat q\rightarrow qt'
\]
given by 
\[
\theta'=g^{-1}:\gamma\rightarrow p\kappa\,,\quad
\phi'=h:\eta\rightarrow q\kappa\,.
\]

The required natural transformation $\psi:1\rightarrow t't$ can be taken
to be the identity. The other one, $\psi':1\rightarrow tt'$ is given on 
the object $(\gamma,g,\kappa,\eta,\eta)$ by the morphism $(g^{-1},1,h)$.
Naturality is assured by virtue of the identities $g'\overline g=g\circ pk'$
and $hh'=qk'\circ\overline h$. The compatibility identities 
(\ref{compatibility}) are easily checked. $\Box$

This completes the verification that $\bB\rightarrow\bC$ is locally an 
equivalence, and hence the proof of Theorem \ref{relate}. $\Box$


\begin{thebibliography}{9}

\bibitem{borceux}
F.~Borceux,
{\em Handbook of Categorical Algebra 1.~ Basic category theory,}
Encyclopedia of Mathematics and its Applications 50, Cambridge University
Press, 1994.

\bibitem{brown}
R.~Brown,
Fibrations of groupoids,
J. Alg. 15 (1970) 103--132.

\bibitem{elmendorf-mandell}
A.~D.~Elmendorf and M.~A.~Mandell, 
Rings, modules, and algebras in infinite loop space theory, 
Adv.~Math.~205 (2006) 163--228.

\bibitem{guillou} 
B.~Guillou,
Strictification of categories weakly enriched in symmetric monoidal 
categories, 
Theory Appl. Categ. 24 (2010) 564--579.

\bibitem{higgins}
P.~J.~Higgins, 
{\em Notes on Categories and Groupoids}, 
van Nostrand Reinhold Mathematical Studies 32, 1971.

\bibitem{lerman}
E.~Lerman,
Orbifolds as stacks?,
Enseign. Math. (2) 56 (2010) 315--363.

\bibitem{leinster}
T.~Leinster,
Basic bicategories,
arXiv:math/9810017v1.

\bibitem{riehl-verity}
E. Riehl and D. Verity,
Kan extensions and the calculus of modules for $\infty$-categories,
arXiv 1507.01460.

\bibitem{schommer-pries}
C.~Schommer-Pries,
Central extensions of smooth 2-groups and a finite-dimensional string
2-group, Geom.~Topol.~15 (2011) 609--676.

\bibitem{tommasini}
M.~Tommasini,
Weak fiber products in bicategories of fractions,
arXiv 412.3295.

\end{thebibliography}
\end{document}